\documentclass[11pt,a4paper]{amsart}

\usepackage[twoside, dvips]{geometry}
\usepackage{amscd}

\newtheorem{prop}{Proposition}
\newtheorem{cor}{Corollary}[prop]
\newtheorem{lemma}{Lemma}

\newcommand{\PP}{\mathbb P}		
\newcommand{\A}{\mathbb A}		
\newcommand{\VV}{\mathbb V}		
\newcommand{\XX}{\mathfrak X}		
\newcommand{\ZC}{\mathfrak S}		
\newcommand{\OO}{\mathcal{O}}		
\newcommand{\Spec}{\mathbf{Spec}}	
\newcommand{\Proj}{\mathbf{Proj}}	

\title{A small resolution for triple covers in algebraic geometry}

\author{Daniele Faenzi}
\email{faenzi@math.unifi.it}
\address{Dipartimento di Matematica ``U.~Dini'', Universit\`a di
  Firenze, Viale Morgagni 67/a, I-50134, Florence, Italy}
\urladdr{http://www.math.unifi.it/\~{ }faenzi/}

\author{Janis Stipins}
\email{jstipins@umich.edu}
\address{Department of Mathematics, University of Michigan, 
  East Hall, 525 East University, Ann Arbor, MI  48109-1109}
\urladdr{http://www.math.lsa.umich.edu/\~{ }jstipins/}

\thanks{2000 {\em Mathematics Subject Classification.} Primary 14E20}

\thanks{The first author was partially supported by the EAGER contract HPRN-CT-2000-00099, and the second author was partially supported by the University of Michigan VIGRE fellowship.  This work was done during the PRAGMATIC 2001 summer school in Catania.}

\title{A small resolution for triple covers in algebraic geometry}

\begin{document}

\begin{abstract}
Given a triple cover $\pi : X \longrightarrow Y$ of varieties,
we produce a new variety $\ZC_X$ and a birational morphism 
$\rho_X : \ZC_X \longrightarrow X$ which is an isomorphism away from
the fat-point ramification locus of $\pi$.  
The variety $\ZC_X$ has a natural interpretation in terms of the data 
describing the triple cover, and the morphism $\rho_X$ has an elegant 
geometric description.
\end{abstract}

\maketitle

\section{Introduction}

The basic fact regarding a triple cover $\pi : X \longrightarrow Y$, proven in 
\cite{miranda-3}, is that any such cover is determined by a rank~2
locally free sheaf $E$ on $Y$ and a global section $\sigma$ of 
$S^3(E)^{*} \otimes \Lambda^2(E)$.  Furthermore, $X$ can be realized
as a subvariety of the geometric vector bundle $\VV(E)$ equipped with
the natural projection to $Y$. 

In this article we give a necessary and sufficient criterion for $X$ to
be realized as a subvariety of a $\PP^1$-bundle equipped with its 
natural projection to $Y$: we show that $X$ can be so realized if and
only if $\pi$ has no fat triple ramification (a fat triple ramification
point of $\pi$ is a point $x \in X$ whose Zariski tangent space in the fibre of
$\pi$ has dimension~2).

Along the way, we show that to any triple cover $\pi : X \longrightarrow Y$, 
one can associate a subvariety $\ZC_X$ of $\PP(E^{*})$ defined in terms
of the global section $\sigma$.  This variety $\ZC_X$ is equipped with 
a birational morphism $\rho_X : \ZC_X \longrightarrow X$, which has a 
nice geometric interpretation.  In fact $\rho_X$ is a sort of small 
resolution: it is the blow-up of a Weil divisor in $X$, and its fibre over any 
fat ramification point of $\pi$ is a $\PP^1$,
but its exceptional set is in general of codimension
larger than 1 in $\ZC_X$.  We construct this resolution first in a
local case, and then we show that the construction globalizes.  Throughout
this article we make extensive use of Miranda's analyses in \cite{miranda-3}.

We note that our main result is not new: a more general statement (with a
correspondingly more technical proof) can be found as Theorem~1.3 in the
beautiful paper \cite{casnati-ekedahl}.  However, it is our hope that our
simple geometric description in the case of triple covers can provide some
insight into the more general case.

The authors wish to thank Ciro Ciliberto, Alberto Calabri, Flaminio Flamini,
Alfio Ragusa, and especially Rick Miranda for their extraordinary efforts
during the PRAGMATIC~2001 summer school in Catania.  The authors also wish
to thank Mike Roth for teaching us about small resolutions, which turned out
to be precisely the right objects for describing the results of our
research.

\section{Some examples of small resolutions}

Let $\A^5$ have coordinates $x,y,z,w,t$, and let $\PP^4$ be its 
projectivization.  Let $X \subset \PP^4$ be the hypersurface 
defined by the equation $xw - yz = 0$.  If we restrict our attention
to the $\PP^3 \subset \PP^4$ where $t=0$, this same equation defines
a smooth quadric $Q \subset \PP^3$; thus $X$ is just the projective 
cone over $Q$ with vertex $[0:0:0:0:1]$.

It is clear that $X$ is smooth away from its vertex.  
If $\epsilon : \widetilde{X} \longrightarrow X$ is the blow-up of the vertex,
then it is easy to see that $\widetilde{X}$ is a smooth variety, and
that the exceptional divisor is isomorphic to $Q = \PP^1 \times \PP^1$.
We are going to describe a method of resolving the singular point
of $X$ with a morphism $\rho : \Gamma \longrightarrow X$, where 
$\Gamma$ is
a smooth variety, but where the exceptional set is a $\PP^1$.
In particular, the exceptional set is ``too small'' to be a divisor;  
thus $\rho$ will be an example of a {\em small resolution}.

We begin by choosing a line $L$ belonging to one of the two rulings of
the quadric $Q$; for ease in computation, we take $L$ to be the line
$x = y = 0$.  Clearly $L$ is a Weil divisor on $Q$.  It is easy to 
check that $L$ is also a Cartier divisor: since there is no point on $Q$
at which $x,y,z,w$ all vanish, the two open sets of $Q$ where $z \neq 0$
and where $w \neq 0$ cover $L$.  On the first set $L$ is defined by
$x=0$, and on the second set $L$ is defined by $y=0$.  

Now let $D \subset X$ be the cone over $L$.  Then $D$ is a Weil divisor
in $X$, but it is not Cartier: $D$ cannot be defined by only one 
equation in any open neighborhood of the origin.
Our small resolution $\rho : \Gamma \longrightarrow X$ is the blow-up of 
$X$ along $D$.  (Note that once we show that $\Gamma$ is not isomorphic to
$X$, we will have proven indirectly that $D$ is not Cartier: the blow-up
of a Cartier divisor is always an isomorphism.)

Since $D$ can be defined by the two equations $x = y = 0$ in $X$, the
blow-up of $X$ along $D$ can be defined as the (closed) graph 
of the rational map $\phi : X --\rightarrow \PP^1$, where 
$\phi ([x:y:z:w:t]) = [x:y]$.  This graph $\Gamma$ is a closed 
subvariety of the
product $X \times \PP^1$, and the morphism 
$\rho : \Gamma \longrightarrow X$ is the restriction
of the first projection map.
Note that there is an open set of $X$ on which
the map $\phi$ agrees with the map $\psi$ 
sending $[x:y:z:w:t]$ to $[z:w]$; that these two maps agree is a 
consequence of the defining equation for $X$.  Also note that at 
any point of $X$ other than its vertex $[0:0:0:0:1]$, at least one
of $\phi$ and $\psi$ is defined.  This observation enables us to
write down the defining equations for $\Gamma \subset X \times \PP^1$:
if $u,v$ are coordinates on the $\PP^1$ factor, then $\Gamma$ is 
defined by $uy - vx = 0$ and $uw - vz = 0$.  It follows that 
$\rho$ is an isomorphism away from the vertex of $X$, and the 
exceptional set over the vertex is $\PP^1$.  It is also easy to 
check that $\Gamma$ is smooth.

We now present an alternative way of describing the variety 
$\Gamma$.  We may regard $\Gamma$ as a subvariety of 
$\PP^4 \times \PP^1$, defined by the three equations 
$xw - yz = 0$, $uy - vx = 0$, and $uw - vz = 0$.  
These three equations may be expressed in a single matrix 
equation:
\[
\Gamma = \left\{ [x:y:z:w:t] \times [u:v] \in \PP^4 \times \PP^1 
	\mbox{ such that }
	\left[ 
		\begin{array}{cc} x & y \\ z & w \end{array} 
	\right]
	\left[
		\begin{array}{c} -v \\ u \end{array}
	\right] = 0
    \right\}.
\]
The first of the three equations is now seen to express the
vanishing of a determinant, 
which is necessary and sufficient for the second and third 
equations to have
a nonzero solution.  
Note that in this description of $\Gamma$, the morphism
$\rho$ is the restriction of the natural projection from 
$\PP^4 \times \PP^1$ to $\PP^4$.

As a second example, we take $X \subset \PP^6$ to be the projective
cone over $\PP^2 \times \PP^1$, with vertex $[\vec{0}:1]$.  By
analogy with the previous example, we consider a divisor 
$D \subset X$ which is the cone over one of the $\PP^1 \times \PP^1$ 
``rulings'' 
of $\PP^2 \times \PP^1$; for example, we can take $D$ to be 
defined by $x_0=x_1=0$.  As before, blowing up this divisor gives
a small resolution $\rho : \Gamma \longrightarrow X$ which is an
isomorphism away from the vertex of $X$, and whose fibre
over the vertex is a $\PP^1$.  A 
computation similar to the previous one gives that $\Gamma$ may
be realized as a subvariety of $\PP^6 \times \PP^1$, using a 
matrix condition:
\[
\Gamma = \left\{ [\vec{x}:t] \times [u:v] \in \PP^6 \times \PP^1 
	\mbox{ such that }
	\left[ 
		\begin{array}{cc} x_0 & x_1 \\ 
				  x_2 & x_3 \\ 
				  x_4 & x_5 \end{array} 
	\right]
	\left[
		\begin{array}{c} -v \\ u \end{array}
	\right] = 0
    \right\}.
\]
This single matrix equation expresses six quadratic conditions:
the vanishing of the three $2 \times 2$ minors, which are the 
defining equations for $X \subset \PP^6$, and the three row 
equations, which come from the blow-up computation.  
As before, the morphism $\rho$ is just the restriction
of the natural projection from $\PP^6 \times \PP^1$ to $\PP^6$.

Note that by simply eliminating the variable $t$ in the 
above matrix description,
we can define $\Gamma$ as a subvariety of $\A^6 \times \PP^1$, 
where $\A^6$ is the finite ($t \neq 0$) part of $\PP^6$.  
This construction is the one that will provide us with a sort of 
universal local picture of our resolution for triple covers.

\section{The local picture of the resolution}

Consider the affine space $\A^4$ with coordinates $A,B,C,D$,
and let $F$ be the free sheaf of rank~2 on this affine space.  Then
$\VV(F)$ is nothing more than the affine space $\A^6$ with 
coordinates $A,B,C,D,z,w$; here $z,w$ are global sections that
generate $F$.  Let $\XX$ be the subvariety of $\VV(F)$ defined by
the three quadrics  
\begin{eqnarray*}
	z^2 & = & Az + Bw + 2(A^2 - BD) \\  
	zw  & = & -Dz - Aw + (BC-AD) \\	    
	w^2 & = & Cz + Dw + 2(D^2 - AC).    
\end{eqnarray*}
By the results in Miranda's paper \cite{miranda-3}, we know that
the projection of $\VV(F)$ to $\A^4$ sending
$(A,B,C,D,z,w)$ to $(A,B,C,D)$ restricts to a triple cover
$\Pi : \XX \longrightarrow \A^4$.

Now, as pointed out in \cite{miranda-3}, the variety $\XX$ is 
determinantal: it is the locus in $\VV(F)$ where 
the matrix
\[
	\left[ 
		\begin{array}{cc} z+A & B \\ 
				  C & w+D \\ 
				  w-2D & z-2A \end{array} 
	\right]
\]
has rank at most one.  By a result in \cite{miranda-3}, the rank of this
matrix is zero if and only if the map $\Pi$ has fat triple
ramification over the point $(A,B,C,D)$; it is clear from the
matrix description that this happens only over the point
$(0,0,0,0)$.

This determinantal representation is familiar: up to a change
of coordinates on $\A^6$, we see that $\XX$ is just the affine cone
over $\PP^2 \times \PP^1$.  Furthermore, we see that the vertex
of this cone -- its only singular point -- is exactly the fat
triple point where $A=B=C=D=z=w=0$.  The temptation to compute
its small resolution is overwhelming, and so we define 
$\Gamma \subset \XX \times \PP^1$
to be the subvariety of $\VV(F) \times \PP^1$
defined by the matrix condition
\[
	\left[ 
		\begin{array}{cc} z+A & B \\ 
				  C & w+D \\ 
				  w-2D & z-2A \end{array} 
	\right]
	\left[
		\begin{array}{c} -v \\ u \end{array}
	\right] = 0.
\]
We know from our previous computation that the natural projection
$\rho : \Gamma \longrightarrow \XX$ is an isomorphism away from
the fat point, and that the fibre over this point
is all of $\PP^1$.  We will refer to the morphism 
$\rho : \Gamma \longrightarrow \XX$ as the {\em resolution of the
triple cover $\Pi$}.

We note that $\Gamma$ comes equipped with a morphism $\phi$ to 
$\A^4 \times \PP^1$: $\phi$ is just the product of 
$\Pi \circ \rho$ with the second projection of $\Gamma$ to $\PP^1$.
We are going to compute the image of $\phi$.  To do 
this, we first note that if $(A,B,C,D,z,w) \times [u:v]$ is a point
of $\Gamma$, then we can solve for $z$ and $w$ in terms of the
other coordinates, using the first two rows of the matrix:
\begin{eqnarray}
\label{eq:phi-inverse}	z &=& B\left(\frac{u}{v}\right) - A \\	    
\nonumber	w &=& C\left(\frac{v}{u}\right) - D.
\end{eqnarray}
Here we assume that both $u$ and $v$ are nonzero; in the case that
either vanishes, the third row of the matrix can be used instead
of one of the first two.  Continuing under the assumption that both
$u$ and $v$ are nonzero, we use the third row of the matrix to
compute that
\[
	-v\left(C\left(\frac{v}{u}\right) - D - 2D\right)
	+u\left(B\left(\frac{u}{v}\right) - A - 2A\right) = 0, 	
\]
and since $uv \neq 0$, we conclude that:
\begin{equation}
	Bu^{3} - 3Au^{2}v + 3Duv^{2} - Cv^{3} = 0.   \label{eq:local-cubic}
\end{equation}
We note that the same equation results from the computations in the
cases where $u=0$ or $v=0$.  

Let $\ZC \subset \A^4 \times \PP^1$ be the subvariety defined by 
equation~(\ref{eq:local-cubic}).
Let $\Pi' : \ZC \longrightarrow \A^4$ be the obvious projection;
this projection is compatible via $\phi$ with the composite map
$\Pi \circ \rho : \Gamma \longrightarrow \A^4$.  In fact, we have
the following result:

\begin{prop}
The morphism $\phi : \Gamma \longrightarrow \ZC$ is an isomorphism
of varieties over $\A^4$.  
\end{prop}
\begin{proof}
The fact that $\Pi \circ \rho = \Pi' \circ \phi$ is clear from
the definition of $\phi$, so we only need to show that $\phi$ is
an isomorphism.  To do this, note that the equations~(\ref{eq:phi-inverse}) 
for $z$ and $w$ in terms
of $A,B,C,D$ define regular functions on all of $\ZC$; this is 
easily checked using equation~(\ref{eq:local-cubic}).  
From the definition
\[
\phi((A,B,C,D,z,w) \times [u:v]) = (A,B,C,D) \times [u:v] ,
\]
we see that $\phi$ is surjective, and also that 
the regular functions for $z$ and $w$ are sufficient to define
the inverse morphism.  Thus $\phi$ is an isomorphism, as needed.
\end{proof}

\begin{cor}
Away from the point $(0,0,0,0) \in \A^4$, the three morphisms 
$\Pi : \XX \longrightarrow \A^4$, 
$\Pi \circ \rho : \Gamma \longrightarrow \A^4$, and 
$\Pi' : \ZC \longrightarrow \A^4$ are isomorphic triple cover maps.
\end{cor}

Now we are in a position to construct a triple cover resolution for 
any sufficiently local triple cover $\pi : X \longrightarrow Y$.  By
``sufficiently local'' we mean that $Y$ is affine, $E$ is a free
sheaf of rank~2 on $Y$, and $X$ is the 
subvariety of $\VV(E)$ defined by the three quadrics
\begin{eqnarray}
\label{eq:local-quadrics}	z^2 & = & az + bw + 2(a^2 - bd) \\  
\nonumber	zw  & = & -dz - aw + (bc-ad) \\	    
\nonumber	w^2 & = & cz + dw + 2(d^2 - ac);    
\end{eqnarray}
here the coefficients $a,b,c,d$ are regular functions on $Y$, and 
$z,w$ are global sections that generate $E$.  It follows from Miranda's 
analysis in \cite{miranda-3} that this is in fact the local situation 
for any triple cover.  

Given such a sufficiently local triple cover, we 
define a morphism $f : Y \longrightarrow \A^4$ by the formula
$f(y) = (a(y), b(y), c(y), d(y))$.  This is equivalent to 
requiring that $f^{*}(A) = a$, and so on; thus we have the 
following commutative diagram:
\[
\begin{CD}
f^{*}\Gamma @>>> \Gamma  \\
@V{f^{*}\rho}VV             @VV{\rho}V \\
X=f^{*}\XX @>>>  \XX                \\
@V{\pi}VV	@VV{\Pi}V         \\
Y     @>f>>    \A^4      
\end{CD}
\]
It is proven in \cite{miranda-3} that the fat points of 
$X \subset \VV(E)$ are precisely the points where $a=b=c=d=0$;
it follows that the morphism $f^{*}\rho$ is an isomorphism away 
from the fat-point ramification locus of $\pi$,
and has a $\PP^1$ fibre over any fat point in $X$.
We will refer to 
$f^{*}\rho$ as the {\em resolution of the triple cover $\pi$}.

In this way we may view the right-hand side of the above diagram
as a sort of universal local picture of our triple cover resolution.  
The reader with some skill in visualizing three-dimensional 
commutative diagrams\footnote{At least, with more skill than the 
authors have in drawing them.} will see that the isomorphism 
$\phi : \Gamma \longrightarrow \ZC$ pulls back via $f$ to an
isomorphism $f^{*}\phi : f^{*}\Gamma \longrightarrow f^{*}\ZC$;
here $f^{*}\ZC$ is the subvariety of $Y \times \PP^1$ defined 
by the the equation
\[
	bu^{3} -3au^2v + 3duv^2 - cv^{3} = 0.
\]
This is in fact a variety over $Y$, whose structure morphism 
$\pi' : f^{*}\ZC \longrightarrow Y$ is none other than $f^{*}\Pi'$.
We can similarly ``pull back'' our other result:

\begin{cor}
Let $B \subset Y$ be the (set-theoretic) image under $\pi$ of the 
fat-point ramification locus in $X$.  Away from $B$, the three
morphisms $\pi : X \longrightarrow Y$, 
$\pi \circ f^{*}\rho : f^{*}\Gamma \longrightarrow Y$, and
$\pi' : f^{*}\ZC \longrightarrow Y$ are isomorphic triple cover maps.
\end{cor}

Thus we reach the following interesting conclusion:

\begin{prop}\label{prop:local-result}
Let $\pi : X \longrightarrow Y$ be a sufficiently local triple cover;
as before, this means that $X$ is defined as a subvariety of a free rank~2
vector bundle on $Y$.  If $\pi$ has no fat-point
ramification, then in fact $X$ is isomorphic as a triple cover to a 
subvariety of a (trivial) $\PP^1$-bundle over $Y$ equipped with the
natural projection.
\end{prop}
\begin{proof}
This is just a restatement of the isomorphism between $X$ and $f^{*}\ZC$
from the previous corollary.
\end{proof}

\section{Geometric description of the resolution}
In this section we are going to describe geometrically the isomorphism
appearing in the preceding proposition.  Along the way, we will also
describe the geometric meaning of the $\PP^1$-bundle appearing there.

To begin, let $\pi : X \longrightarrow Y$ be any triple cover.  
For convenience of notation, we set $F = \pi_*(\OO_X)$.  Following 
\cite{miranda-3}, we have that 
$F = \OO_Y \oplus E$,
where $E$ is a rank~2 locally free sheaf on $Y$.  If $U \subset Y$
is any open set over which $E$ is generated freely by two sections
$z,w$, then over $U$ we have that $X$ is defined as a subvariety
of $\VV(E)$ by the quadrics~(\ref{eq:local-quadrics}); this makes 
sense, because local sections of $E$ correspond to local coordinates
on $\VV(E)$.  Over a fixed $y \in Y$, the fibre of $\VV(E)$ is an 
affine plane, and the quadrics~(\ref{eq:local-quadrics}) cut out one, 
two, or three points in this plane.  

Our idea is to consider the function on $X$ that
is defined by sending a point in a fibre of $\pi$ to the line through the 
other two points in the fibre; clearly we need to work a bit to understand 
this idea.  For one thing, this definition only makes sense for
fibres containing three distinct points of $X$.  Still, we may hope
to define a rational map on $X$ whose locus of indeterminacy is 
contained in the ramification locus of $\pi$.  The greater difficulty
is understanding what the range of this function should be: we need to map
to a bundle whose fibre over a fixed point in $Y$ is the $\PP^1$ of
lines in the fibre of $\VV(E)$.

It turns out that we can make this idea work by considering the 
inclusion $E \hookrightarrow \OO_Y \oplus E = F$.  This inclusion
allows us to identify the fibre
of $\VV(E)$ with the ``finite part'' of the fibre of $\PP(F)$, 
which is a projective plane.  The ``line at infinity'' in the fibre 
of $\PP(F)$ is identified with the fibre of $\PP(E)$.  (These
identifications follow from applying the functors $\Spec$ and 
$\Proj$ to the stated inclusion.)  As sets, we have that 
$\PP(F) = \VV(E) \cup \PP(E)$, and so we may view $X \subset \VV(E)$
as a subvariety of $\PP(F)$ that does not intersect $\PP(E)$.

Now we are in a position to describe our putative rational map on $X$.
Over a point $y$ not in the branch locus of $\pi$, the fibre of $\pi$
consists of three distinct points $x_1,x_2,x_3$.  The line through
any two of these points, say $x_2$ and $x_3$, is a line in the fibre
of $\PP(F)$ over $y$.  Such a line corresponds to a point $p_{2,3}$ in the
fibre of $\PP(F^{*})$ over $y$.  There is a natural projection 
$\PP(F^{*}) --\rightarrow \PP(E^{*})$ which dualizes the inclusion 
$\PP(E) \hookrightarrow \PP(F)$; over the point $y$, this map is 
the projection whose center is the point corresponding to the line
at infinity in the fibre of $\PP(F)$ over $y$.  We have that $p_{2,3}$ is
never equal to the center of this projection, because $X$ does not meet
$\PP(E)$; thus we can project $p_{2,3}$ to a point $q_{2,3}$ in 
the fibre of $\PP(E^{*})$ over $y$.  We define a rational map 
$\psi : X --\rightarrow \PP(E^{*})$ by setting $\psi(x_1) = q_{2,3}$, 
and similarly for $x_2,x_3$.

We have the following result:
\begin{prop}
If $\pi : X \longrightarrow Y$ is sufficiently local, then the map
$\psi$ is in fact rational, and its image is contained in 
$f^{*}\ZC \subset \PP(E^{*}) = Y \times \PP^1$.  The restricted map
\[
	\psi : X --\rightarrow f^{*}\ZC
\]
is birational, and the
isomorphism $f^{*}\phi : f^{*}\Gamma \longrightarrow f^{*}\ZC$ is
the resolution of indeterminacy of this birational map.
\end{prop}
All of these claims can be checked easily (by the reader!) 
once we establish the 
expression for $\psi$ in terms of local coordinates on 
$X \subset \VV(E)$:

\begin{lemma}
The local expression for $\psi$ over a point $y \in Y$ is
\begin{eqnarray*}
	\psi(y \times (z,w)) &=& y \times [z + a(y) : b(y)] \\
		  &=& y \times [c(y) : w + d(y)] \\
		  &=& y \times [w-2d(y) : z-2a(y)],
\end{eqnarray*}
where $z,w$ are coordinates on the fibre of $\pi$ over $y$, and $a,b,c,d$
are the sections of $\OO_Y$ appearing as coefficients in the 
equations~(\ref{eq:local-quadrics}).
\end{lemma}
Note that
the equivalence of the three expressions for $\psi$ is a consequence
of the equations~(\ref{eq:local-quadrics}) which define $X$ as a 
subvariety of $\VV(E)$.
\begin{proof}
In order to proceed, we need to recall a fact from \cite{miranda-3}
regarding the sheaf $E$: the $\OO_Y$-algebra $\pi_{*}(\OO_X)$ is in 
fact a rank~3 $\OO_Y$-module, and $E$ is the rank~2 submodule 
consisting of sections that have zero trace over $\OO_Y$.  This means
that if we take local generators $z,w$ of $E$ as local coordinates
on $\VV(E)$, then the vector sum of the points $(z_i,w_i)$ 
in any fibre of $\pi$ must be zero.  

Now we can prove the lemma.  Let $y \in Y$ be any point not in the
branch locus of $\pi$.  Then the fibre of $\pi$ over $y$ consists of
three distinct points, which we denote $(z_1,w_1),(z_2,w_2),(z_3,w_3)$.
The fibre of $\VV(E)$ over $y$ is an affine plane with coordinates
$z,w$; inside this plane, the line containing $(z_2,w_2)$ and $(z_3,w_3)$ is
given by
\[
-(w_3 - w_2)z + (z_3 - z_2) w + (z_2 w_3 - z_3 w_2) = 0. 
\]
In local coordinates, then, we have
\[
\psi(z_1 , w_1) = [ -(w_3 - w_2) : (z_3 - z_2) ];
\]
this is understood to be a point in the fibre of $\PP(E^{*})$ over $y$.

We claim that this expression agrees with the first one given in the 
statement of the lemma.  To see this, we use the 
equations~(\ref{eq:local-quadrics}) and the zero trace observation above
to compute that
\begin{eqnarray*}
	\left( z_1 + a(y) \right) (z_3 - z_2) &=& z_1 z_3 - z_1 z_2 + a(y)(z_3 - z_2) \\
	&=& z_1 z_3 - z_1 z_2 + \left( z_{3}^2 - z_{2}^2 - b(y)(w_3 - w_2) \right) \\
	&=& (z_1 + z_2 + z_3)(z_3 - z_2) - b(y)(w_3 - w_2) \\
	&=& 0 - b(y)(w_3 - w_2) .
\end{eqnarray*}
This proves that $\psi(z_1,z_2) = [ z_1 + a(y) : b(y) ]$ on the open set where
$\pi$ is unramified and where this expression is defined.  
Since we only require $\psi$ to be a rational map, the lemma is proved.
\end{proof}

This result shows that the locus of indeterminacy of $\psi$
is precisely the fat-point ramification locus of $\pi$, which in general
is a proper subset of the ramification locus of $\pi$.  This is 
consistent with the fact that reasonable definitions of the rational 
map $\psi$ can be made for double ramification points and for 
curvilinear triple ramification points; in these cases the Zariski
tangent spaces to the ramification points determine lines in the
fibres of $\pi$.  At a fat point, the dimension of the Zariski tangent
space in the fibre is equal to 2, so there is no reasonable way
to define $\psi$ at such a point.

\section{Globalization}

Now we are going to define our resolution for an arbitrary triple cover 
$\pi : X \longrightarrow Y$.  The idea is straightforward: we know
from \cite{miranda-3} that $Y$ is covered by open affine sets $Y_i$ for
which the restricted triple covers $\pi : X_i \longrightarrow Y_i$
are sufficiently local, and we
have already defined the resolution for sufficiently local triple
covers.  It remains to check that these local definitions patch
together compatibly to define a global resolution.  

Recall that in defining
the universal local resolution $\rho : \Gamma \longrightarrow \XX$,
we constructed $\Gamma$ as a subvariety of $\VV(F) \times \PP^1$.
The unidentified factor of $\PP^1$ is an obstruction to globalization:
if each sufficiently local resolution variety is defined as a 
subvariety of $\VV(E)_{|Y_i} \times \PP^1$, then it is not clear how to interpret
the second factor as the restriction of a globally defined object.
Fortunately, we have seen how to remedy this: using the isomorphism
$\phi : \Gamma \longrightarrow \ZC$, we will take $\ZC$ to be our
resolution variety instead of $\Gamma$.  Then we take the resolution
variety for a sufficiently local triple cover to be $\ZC_i = f_{i}^{*}\ZC$ 
instead of $f_{i}^{*}\Gamma$.  In the previous section we showed that
each variety $\ZC$ is naturally a subvariety of 
$\PP(E^{*})_{|Y_i}$.  Thus we may hope to patch together the varieties
$\ZC_i$ to construct a subvariety $\ZC_X$ of $\PP(E^{*})$.

Now we will invoke a beautiful result of Miranda from \cite{miranda-3} 
to finish our construction.  Miranda shows that any triple cover 
$\pi : X \longrightarrow Y$ is determined by a rank~2 locally free sheaf 
$E$ on $Y$ and a global section $\sigma$ of $S^3(E)^{*} \otimes \Lambda^2(E)$.
In fact, Miranda shows that if $a,b,c,d$ are the coefficients appearing
in the quadrics~(\ref{eq:local-quadrics}) that define $X_i$ as a 
subvariety of $\VV(E)_{|Y_i}$, then the local expression for $\sigma$ over
$Y_i$ is
\[
	-b(z^3)^{*} + a(z^2)^{*}w^{*} - dz^{*}(w^2)^{*} + c(w^3)^{*}.
\]
Using the natural isomorphism $S^3(E)^{*} \cong S^3(E^{*})$, we get the
following local expression for $\sigma$:
\[
	-\frac{1}{6}b(z^{*})^3 + \frac{1}{2}a(z^{*})^{2}w^{*} - \frac{1}{2}dz^{*}(w^{*})^{2} + \frac{1}{6}c(w^{*})^{3}.
\]
Up to a constant factor, this is just the cubic defining $\ZC_i$ as a 
subvariety of $\PP(E^{*})_{|Y_i}$.  Since $\sigma$ is a global section,
we conclude that the varieties $\ZC_i$ must patch together to form a
variety $\ZC_X \subset \PP(E^{*})$.  It is clear that the structure 
morphisms patch together compatibly, and so we have the following result:

\begin{prop}
Let $\pi : X \longrightarrow Y$ be any triple cover, and
let $E$ be a rank~2 locally free sheaf on $Y$ such that $X \subset \VV(E)$. 
Then there is a variety $\ZC_X \subset \PP(E^{*})$
and a birational morphism $\rho_X : \ZC_X \longrightarrow X$
which is an isomorphism away from the fat-point ramification locus of
$\pi$, and whose fibre over every fat point is a $\PP^1$.
\end{prop}

We refer to the morphism $\rho_X : \ZC_X \longrightarrow X$ as the {\em 
resolution of the triple cover $\pi$}.  Now we get the following 
global result:

\begin{prop}
Let $\pi : X \longrightarrow Y$ be a triple cover.  $X$ is isomorphic
as a triple cover to a subvariety of a $\PP^1$-bundle on $Y$ equipped
with the natural projection if and
only if $\pi$ has no fat-point ramification.
\end{prop}
\begin{proof}
One implication is the global version of Proposition~\ref{prop:local-result};
the other implication
follows from the fact that the fibre of a subvariety of a $\PP^1$-bundle
cannot have a two-dimensional Zariski tangent space.
\end{proof}


\bigskip

\end{document}